\documentclass [final]{siamltex}
\usepackage{amsfonts,amssymb}

\newcommand{\commentout}[1]{}
\newcommand{\ba}{\begin{array}}
        \newcommand{\ea}{\end{array}}
\newcommand{\bdm}{\begin{displaymath}}
        \newcommand{\edm}{\end{displaymath}}
\newcommand{\ben}{\begin{enumerate}}
        \newcommand{\een}{\end{enumerate}}
\newcommand{\beq}{\begin{equation}}
        \newcommand{\eeq}{\end{equation}}
\newcommand{\bfg} {\begin{figure}[p]}
        \newcommand{\efg} {\end{figure}}
\newcommand{\bi} {\begin {itemize}}
        \newcommand{\ei} {\end {itemize}}
\newcommand{\bqn}{\begin{eqnarray}}
        \newcommand{\eqn}{\end{eqnarray}}
\newcommand{\bqs}{\begin{eqnarray*}}
        \newcommand{\eqs}{\end{eqnarray*}}
\newcommand{\bsl} {\begin{slide}[8.8in,6.7in]}
        \newcommand{\esl} {\end{slide}}
\newcommand{\bss} {\begin{slide*}[9.3in,6.7in]}
        \newcommand{\ess} {\end{slide*}}
\newcommand{\btb} {\begin {table}}
        \newcommand{\etb} {\end {table}}
\newcommand{\bvb}{\begin{verbatim}}
        \newcommand{\evb}{\end{verbatim}}

\newcommand{\m}{\mbox}
\newcommand{\mat}[1]{{{\left[ \ba #1 \ea \right]}}}
\newcommand{\cas}[1]{{{\left \{ \ba #1 \ea \right. }}}

\newcommand{\refe}[1] {{{(\ref {#1})}}}

\def\pmb#1{\setbox0=\hbox{$#1$}%
   \kern-.025em\copy0\kern-\wd0
   \kern.05em\copy0\kern-\wd0
   \kern-.025em\raise.0433em\box0 }

\def\f{{\mathfrak F}}
\def\gamma{{h}}
\def\p{{\mathfrak P}}
\def\phi{{P}}
\def\alpha{{Q}}
\renewcommand{\tilde}[1] {{\bar #1}}

\def\phitwo     {{\p^{[2]}(x)}}

\def\L      {{\mathbf L}}
\def\X      {{\textbf X}}
\def\bruna  {{\left [\begin {array} {ccccc}
        0&1&&&\\
        &0&1&&\\
        & &0&\ddots&\\
        &&&\ddots&1\\
        &&&&0
        \end {array} \right]}}
\def\brunb  {{\left[ \begin {array}{c}
        0\\
        0\\
        \vdots\\
        0\\
        1
        \end {array} \right]}}

\newcommand{\matphi} [2]{{\left [ \begin {array} {c}
        [\phi_1]_{(#2)} \\\hline
        [\phi_2]_{(#2)} \\\hline
        \vdots \\\hline
        [\phi_#1]_{(#2)}
        \end {array} \right ]}}

\newcommand{\matx}[2]   {{\left [ \begin {array} {c}
    [\L^0 \phi]_{(#2)} \\\hline
    [\L^1 \phi]_{(#2)} \\\hline
    \vdots \\\hline
    [\L^{#1-1} \phi]_{(#2)}
    \end {array} \right ]}}

\newlength {\wid} \newlength {\hei}

\newtheorem{algorithm}[theorem]{Algorithm}
\newtheorem{example}{Example}

\def\eop{{\hfill $\blacksquare$}}

\title {A method for computing quadratic Brunovsky forms}
\author {Wen-Long Jin\thanks{Institute of Transportation Studies, University of California, 522 Social Science Tower, Irvine, CA 92697-3600, USA. Email: wjin@uci.edu}}

\begin {document}
\maketitle

\begin{abstract}
In this paper, for continuous, linearly-controllable quadratic control systems with a single input, we propose an explicit, constructive method for studying their Brunovsky forms, initially studied in \cite{kang1992}. In this approach, the computation of Brunovsky forms and transformation matrices and the proof of their existence and uniqueness are carried out simultaneously. In addition, we show that quadratic transformations in \cite{kang1992} can be simplified to prevent multiplicity in Brunovsky forms. We further extend our method for studying discrete quadratic systems. Finally, computation algorithms for both continuous and discrete systems are summarized, and examples demonstrated.
\end{abstract}



\pagestyle{myheadings}
\thispagestyle{plain}
\markboth{W.-L. JIN}{QUADRATIC BRUNOVSKY FORMS}

\section {Introduction}
Linear control systems can be continuous in time $t$,
\bqn
\dot \xi=A \xi+b\mu, \label{lin_con}
\eqn
or discrete
\bqn
\xi(t+1)=A \xi(t)+b\mu(t),\label{lin_dis}
\eqn
where coefficients $A\in\Re^{n\times n}$ and $b\in\Re^{n\times1}$ are generally constant, and state variable, $\xi$($\in\Re^n$), and control variable, $\mu$($\in \Re$), are continuous or discrete respectively. When  controllable,  both \refe{lin_con} and \refe{lin_dis} admit the Brunovsky form \cite{brunovsky1970}, which is derived from controller form, under the following linear change of coordinates and state feedback
\bqn
\xi=T x,\qquad \mu=u+x^T v, \label{lin_tran}
\eqn
where $x$ and $u$ are new state and control variables respectively, and $T\in\Re^{n\times n}$ and $v\in\Re^{n\times1}$ are the transformation matrices (vectors) (refer to Chapter 3 of \cite{kailath1980}). In the linear Brunovsky form of \refe{lin_con} and \refe{lin_dis}, $A$ and $b$ have the following forms:
\bqn
A=\bruna,\quad b=\brunb. \label{def:AB}
\eqn

In an attempt to extending the Brunovsky forms for non-linear systems, which generally do not even have controller forms, Kang and Krener \cite{kang1992} studied continuous linearly controllable quadratic control systems with a single input, which can be written as (We use different notations from \cite{kang1992} for the purpose of clearly presenting our method of computation.)
\bqn
\dot \xi=A \xi+b\mu+\f^{[2]}(\xi)+G \xi \mu+ O (\xi,\mu)^3,\label{sys:consin}
\eqn
where $\f^{[2]}(\xi)=(\xi^T F_1 \xi,\cdots,\xi^T F_n \xi)^T$ is a vector of $n$ quadratic terms with symmetric $n\times n$ matrices $F_i$'s ($i=1,\cdots,n$), $G\xi\mu=(G_1 \xi \mu, \cdots,  G_n\xi \mu)^T$ is a vector of $n$ bilinear terms with $G\in \Re^{n\times n}$, and $O (\xi,\mu)^3$ includes all terms $\xi^a\mu^b$ with $a+b\geq 3$. Moreover, the following additional assumptions are made for this system. First, the coefficients in \refe{sys:consin} and \refe{sys:dissin} $A$, $b$, $F_i$ ($i=1,\cdots, n$), $G$ and $h$ are assumed to be time-invariant. Second, the two systems are assumed to be linearly controllable; by linearly controllable we mean that the linear parts of the two quadratic systems are controllable and, as a result, the linear parts can be transformed into the Brunovsky form with \refe{lin_tran}.

Kang and Krener \cite{kang1992} first defined quadratically state feedback equivalence up to second-order, or quadratic equivalence for short, under the following quadratic change of coordinates and state feedback,
\bqn
\xi=x+\phitwo+ O (x)^3, \qquad \nu=\mu+x^T Q x+r x \mu+ O (x,\mu)^3, \label{qu2_tran}
\eqn
in which $\phitwo=(x^T P_1 x,\cdots,x^T P_n x)^T$ is a vector of $n$ quadratic terms, and transformation matrices include symmetric $P_i\in \Re^{n\times n}$ ($i=1,\cdots,n$), symmetric $Q\in \Re^{n\times n}$, and $r\in\Re^{1\times n}$. We can see that \refe{qu2_tran} are equivalent to
\bqn
\xi=x+\phitwo+ O (x)^3, \qquad \mu=\nu-x^T Q x-r x \nu+ O (x,\nu)^3, \label{qu_tran}
\eqn
and hereafter refer to transformations of \cite{kang1992} as \refe{qu_tran}.
Then, from all quadratically equivalent systems of a general system \refe{sys:consin}, two types of Brunovsky forms were defined in \cite{kang1992}. In type I forms, the nonlinear terms are reduced into a number of quadratic terms $x_i^2$; i.e., there are no the cross terms in $x_i$ and $x_j$ ($i\neq j$) or bilinear terms, $x\nu$. In type II forms, only bilinear terms are kept. In both types of Brunovsky forms, the number of non-zero nonlinear terms is $n(n-1)/2$, compared to $n^2(n+3)/2$ for a general quadratic system.

In this paper, we propose a new method for studying the Brunovsky forms, first for \refe{sys:consin} under transformations \refe{qu_tran}. This method can be carried out in three steps: first, we find the relationships between the coefficients of \refe{sys:consin}, the coefficients of its quadratically equivalent systems, and the corresponding transformation matrices; second, from these relationships, we derive a mapping from the coefficients of \refe{sys:consin} and its equivalent systems to the transformation matrix $P_1$, which can be considered as the necessary condition that all equivalent systems of \refe{sys:consin} should satisfy; third, we show how to compute, from the necessary condition, Brunovsky forms as well as the corresponding transformations. With this method, we find that \refe{sys:consin} admits the same two types of Brunovsky forms defined in \cite{kang1992}. In contrast, our approach, constructive in nature, is capable of computing the Brunovsky forms and the transformation matrices simultaneously.

We further show that quadratic transformations, \refe{qu_tran}, can be simplified by setting $r=0$ to
\bqn
\xi=x+\phitwo+ O (x)^3, \qquad \mu=\nu-x^T \alpha x+ O (x,\nu)^3. \label{qud_tran}
\eqn
Still defining the quadratic equivalence in the sense of \cite{kang1992}, the new transformations prevent multiple solutions of type I or type II Brunovsky forms. I.e.,
the Brunovsky form of each type and the corresponding transformation matrices $P_i$ ($i=1,\cdots,n$) and $Q$ are uniquely determined by the original system. Moreover, we apply the same method, but with the simplified transformations defined in \refe{qud_tran}, and study the following discrete system
\bqn
\xi(t+1)=A\xi(t)+b\mu(t)+\f^{[2]}(\xi(t))+G\xi(t)\mu(t)+\gamma \mu^2(t)+O(\xi,\mu)^3,\label {sys:dissin}
\eqn
where, similarly, $\f^{[2]}(\xi(t))$ and $G\xi(t)\mu(t)$ are a vector of quadratic terms and a vector of bilinear terms respectively. \footnote{Note that the discrete system \refe{sys:dissin} contains a term quadratic in the control variable, $\gamma \mu^2(t)$, where $\gamma\in\Re^{n\times 1}$, and \refe{sys:dissin} has $n^2(n+3)/2+n$ non-zero nonlinear terms.} We find that for \refe{sys:dissin} there exists only one type of Brunovsky form consisting of $n(n+1)/2$ bilinear terms, which corresponds to type II Brunovsky forms of continuous systems.

The rest of this paper is organized as follows. After reviewing the Brunovsky form of linear systems (Section 2), we study the Brunovsky forms of continuous quadratic systems in Section 3 and of discrete quadratic systems in Section 4. In Section 5, we summarize our method into two computation algorithms: one for continuous systems, and the other for discrete systems. We conclude our paper in Section 6.

\section {Review: computation of the Brunovsky form of a controllable linear control system}

We first review the computation of the Brunovsky form of a continuous controllable linear control system \refe{lin_con}. \footnote{For further explanation, refer to \cite{brunovsky1970} and Chapter 3 of \cite{kailath1980}.} The procedure and the results also apply to the discrete system \refe{lin_dis}. Computation of the Brunovsky form for \refe{lin_con} comprises of two steps: first, the linear system is transformed into the controller form with a linear change of coordinates given by the first equation in \refe{lin_tran}; second, the controller form is further reduced into the Brunovsky form with a state feedback given by the second equation in \refe{lin_tran}.

The control matrix for \refe{lin_con} is defined as
\bqn
C&=&[A^{n-1}b \quad \cdots \quad Ab\quad  b].
\eqn
Since \refe{lin_con} is controllable, rank($C$) is $n$ and, therefore, $C$ is invertible.
Denote $d$ as the first row of $C^{-1}$, then we can compute the transformation matrix $T$ as
\bqn
T&=&\mat{{c} d \\d A \\\vdots\\d A^{n-2}\\d A^{n-1}},
\eqn
and with the linear change of coordinates $\xi=T x$, \refe{lin_con} can be transformed into the following controller form,
\bqn
\dot x&=&\bar A x +\bar b \mu,\label{controller}
\eqn
in which
\bqn
\bar A=T^{-1} A T=
\left [\begin {array} {ccccc}
        0&1&&&\\
        0&0&1&&\\
        \vdots& &0&\ddots&\\
        \vdots&&&\ddots&1\\
        -v_1&-v_2&\cdots&\cdots&-v_n
        \end {array} \right], \qquad
\bar b&=& T^{-1} b=\brunb.
\eqn

Then, using a linear state feedback $\mu=u+x^Tv$, in which,
\bqn
v&=&\mat{{c}v_1\\v_2\\\vdots\\v_{n-1}\\v_n},
\eqn
we can transform the controller form \refe{controller} into the following Brunovsky form,
\bqn
\dot x&=&\bruna x+\brunb u.
\eqn

From the derivation above, we can see that the Brunovsky form exists and is unique, and the transformations can be explicitly computed. In addition, this procedure can be easily integrated in applications related to the Brunovsky form. In the same spirit, we propose a method for directly computing quadratic Brunovsky forms and the corresponding transformations in the following sections.

\section {Computation of continuous quadratic Brunovsky forms}

Note that, with linear transformations in \refe{lin_tran}, quadratic systems \refe{sys:consin} and \refe{sys:dissin} will not change their forms except their coefficients. Without loss of generality, therefore, we hereafter assume that $A$ and $b$ in \refe{sys:consin} ( and also in \refe{sys:dissin}) have already been transformed into the forms defined by \refe{def:AB}.
In this section, we study Brunovsky forms of continuous system \refe{sys:consin} under transformations \refe{qu_tran}, in a constructive manner. First, we find relationships between coefficients of \refe{sys:consin}, coefficients of its quadratically equivalent systems, and the corresponding transformation matrices. Second, from these relationships, we derive a mapping from the coefficients of \refe{sys:consin} and those of the equivalent systems to transformation matrix $P_1$; this mapping is a necessary condition that all the equivalent systems should satisfy. Third, we show how to obtain two types of Brunovsky forms as well as the corresponding transformations.

\subsection{Quadratically equivalent system of \refe{sys:consin}}

\begin{definition}\label{def1}
Assume $A\in\Re^{n\times n}$ is in the linear Brunovsky form given by \refe{def:AB}, we define a linear operator $\L:\Re^{n\times n}\to\Re^{n\times n}$ as
\bqn
\ba{lcl}
\L^0 P  &=& P, \\
\L P&=&A^T P+PA,\\
\L^{i+1} P &=&\L \cdot \L^i P, \quad i=0,1,\cdots.
\ea\label{equ2:6}
\eqn
\end{definition}
\begin{description}
\item[Properties of $\L$.] The linear operator $\L$ has the following properties:
\ben
\item $\L^i \phi=0$ when $i\geq 2n-1$.
\item The nullity of $\L$ is $n$; if $P \in ker(\L)$; i.e., $\L P=0$, then $P$ can be written as
\bqs
P_{ij}&=&\cas{{ll}0,&i+j\leq n,\\(-1)^ip_{i+j-n},&\m{otherwise,}}
\eqs
in which $p_1,\cdots,p_n$ are independent.
\item $\L$ is not invertible.
\item If $P$ is symmetric, $\L P$ is symmetric.

\een
\end{description}

\begin{theorem}\label{theorem:con-equiv}
The continuous quadratic system \refe{sys:consin} is equivalent, in the sense of \cite{kang1992}, under the quadratic transformations given by \refe{qu_tran}, to a quadratic system, whose $i$th ($i=1,\cdots,n$) equation is
\bqn
\dot x_i &=& x_{i+1} +b_i \nu +x^T \tilde F_i x +\tilde G_i x \nu + O  (x,\nu)^3, \label {equ:normalform}
\eqn
where $x_{n+1}(t)\equiv 0$ is a dummy state variable, $\bar F_1,\cdots,\bar F_n$  are symmetric, $\bar G_i$ is the $i$th row of the matrix $\bar G$, and the coefficients $\bar F_i$ ($i=1,\cdots,n$) and $\bar G$ are defined by
\bqn
\tilde F_i&=&F_i+\phi_{i+1}-\L \phi_i-b_i\alpha, \label{f:iter}\\
\tilde G_i     &=& G_i-2b^T\phi_i-b_i r,\label {equ:second}
\eqn
where $\phi_{n+1}\in \Re^{n\times n}$ is a zero dummy transformation matrix, and $b_i$ is the $i$th element of $b$ defined in \refe{def:AB}.
\end{theorem}

{\em Proof.}
Plugging the transformations defined in \refe{qu_tran} into \refe{sys:consin}, we obtain
\begin {eqnarray*}
\dot x+\frac d {dt} \phitwo&=&Ax+b \nu+\f^{[2]}(x)+G x \nu \\
    &&+A\phitwo -bx^TQx-brx\nu+ O (x,\nu)^3,
\end {eqnarray*}
of which the $i$th ($i=1,\cdots,n$) differential equation can be written as
\bqs
\dot x_i+\frac d {dt} (x^T \phi_i x) &=&x_{i+1}+b_i \nu+x^TF_ix+G_ix \nu\\&&+ x^T\phi_{i+1}x -b_i x^T\alpha x -b_irx\nu
    + O (x,\nu)^3.
\eqs
After plugging $\dot x= Ax+b\nu + O (x,\nu)^2$ into the term $\frac d {dt} (x^T \phi_i x)$ and collecting all terms whose orders are higher than two, we then obtain the equivalent system \refe{equ:normalform} with the coefficients given in \refe{f:iter} and \refe{equ:second}. \eop

\begin{description}
\item[Remarks] Equations \ref{f:iter} and \ref{equ:second} present the relationships between the coefficients of \refe{sys:consin}, those of its quadratically equivalent system \refe{equ:normalform}, and the corresponding transformation matrices. We further simplify these relationships in the following subsections.
\end{description}

\subsection{A necessary condition for quadratically equivalent systems}

\begin{lemma}\label{lemma:pg}
The mapping from the coefficients, $G$ and $\bar G$, and $r$ to the transformation matrices $P_i$ ($i=1,\cdots,n$) is given by
\begin {eqnarray}
\matphi{n}{n}&=&\frac 12 G-\frac 12 \tilde G -\frac 12 b r,\label {equ:fourth}
\end {eqnarray}
where the operator $[\cdot]_{(n)}$ takes the $n$th row from its object.
\end{lemma}

{\em Proof.} Since $b_i=0$ ($i=1,\cdots, n-1$) and $b_n=1$, we obtain, from (\ref {equ:second}),
\bqs
[\phi_i]_{(n)} &=& b^T \phi_i = \frac 1 2 G_i-\frac 12 \tilde G_i -\frac 12 b_i r .
\eqs
Thus we have \refe{equ:fourth}. \eop

\begin{lemma}\label{lemma:pf}
The mapping from the coefficients $F_i$ and $\bar F_i$ ($i=1,\cdots,n$) to the transformation matrices $P_i$ ($i=1,\cdots,n$) and $Q$ is given by
\bqn
\phi_{i+1} &=&\L^i\phi_1-\sum_{j=0}^{i-1} \L^j F_{i-j} +\sum_{j=0}^{i-1}\L^j \tilde F_{i-j}, \quad i=1,\cdots, n-1, \label {equ:third} \\
\alpha  &=&\sum_{j=0}^{n-1} \L^j (F_{n-j}-\tilde F_{n-j})-L^n \phi_1. \label{equ:alpha}
\eqn
\end{lemma}
{\em Proof.} When $i=1,\cdots,n-1$, from \refe{f:iter}, we have
\bqs
\phi_{i+1} &=&\L \phi_i-F_i+\tilde F_i.
\eqs
After iterating this equation with respect to $i$, we obtain \refe{equ:third}. When $i=n$,  \refe{f:iter} can be written as
\bqs
\alpha  &=&F_n-\tilde F_n -\L \phi_n.
\eqs
Since $P_n$ is given in \refe{equ:third}, we then have \refe{equ:alpha}. \eop

\begin{definition}\label{def:xc}
With $A\in\Re^{n\times n}$ and $\L$ defined in \refe{def:AB} and Definition \ref{def1} respectively, we define a series of linear operators $\X_i:\Re^{n\times n}\to \Re^{n\times n}$ ($i=0,1,\cdots$) by
\begin {eqnarray}
\X_0 \phi &=&\matx nn, \quad \X_i \phi = (A^T)^i \X_0 \phi. \label{def:Xi}
\end {eqnarray}
\end{definition}
\begin{description}
\item[Properties of $\X_i$.] The linear operators $\X_i$ ($i=0, 1,\cdots$) have the following properties:
\ben
\item $\X_0$ transforms a diagonal matrix into a skew-triangular matrix of the following structure. For a $k$th ($k=-n+1,\cdots,n-1$) diagonal matrix $P$, we denote the $k$th diagonal elements by $p_l=P_{(|k|-k)/2+l,(|k|+k)/2+l}$ ($l=1,\cdots,n-|k|$). Then $\X_0 P$ becomes a diagonal matrix with the following properties: $(\X_0 P)_{n-(|k|-k)/2-l+1,(|k|+k)/2+l}=p_l$, all $(\X_0 P)_{i,j}$ for $i+j=n+k+1+2m$, where $m\geq 1$ and $n+k+1+2m\leq 2n$\footnote{Note that, for $n=2$, there is no solution to $m$. In this case, all entries are zeros except $(\X_0 P)_{ij}$ for $i+j=n+k+1$.}, are determined by $p_l$, and the other entries are zeros.
\item From the preceding property, $\X_0$ transforms a lower-triangular matrix into a full matrix and transforms an upper-triangular matrix into a lower skew-triangular matrix $\Delta$, defined by
\bqn
 \Delta_{ij}=0, \qquad \m{ when $i+j\leq n+1$} \label{def:Delta}.
\eqn
\item From Properties 1 and 2, we have that the nullity of $\X_0$ is 0. Therefore, $\X_0$ is invertible.
\item From the definition of $\L$, $\X_i=0$ when $i\geq n$.
\item From the definition of $\X_i$ ($i=1,\cdots,n-1$), $\X_i P$  can be obtained by shifting $\X_0 P$ down by $i$ rows and replacing the first $i$ rows by zeros. From Property 1, therefore, $\X_i$ transforms a diagonal matrix $P$, whose main diagonal elements are $p_1,\cdots,p_n$ into a lower skew-triangular matrix of the following structure: $(\X_i P)_{n-l+i+1,l}=p_{l}$ for $l=i+1,\cdots,n$, and all other elements except $(\X_i P)_{i,j}$ for $i+j=n+i+2m$, where integer $m\geq 1$ and $n+k+1+2m\leq 2n$, are zeros.
\een
\end{description}

\begin{theorem} \label{con:p1}
The mapping from $F_i$, $G$, $\bar F_i$, $\bar G$ ($i=1,\cdots,n-1$), and $r$ to $P_1$ is given by
\begin {eqnarray}
\phi_1&=&\X_0^{-1}\left (\sum _{i=1} ^{n-1} \X_i F_i+\frac 12 G -\sum _{i=1} ^{n-1} \X_i \tilde F_i - \frac 12 \tilde G-\frac 12 b r \right). \label {equ:important}
\end {eqnarray}
\end{theorem}
{\em Proof.}
From (\ref {equ:third}) we can find the $n$th row of $\phi_{i+1}$ ($i=1,\cdots, n-1$) as
\bqs
[\phi_{i+1}]_{(n)} &=& [\L ^i \phi_1]_{(n)} -\sum _{j=0}^ {i-1}[\L^j F_{i-j}]_{(n)} +\sum _{j=0}^ {i-1}[\L^j \tilde F_{i-j}]_{(n)}.
\eqs
Hence, we have
\begin {eqnarray}
\matphi{n}{n}&=&\X_0 \phi_1 - \sum _{i=1} ^{n-1} \X_i F_i +\sum _{i=1} ^{n-1} \X_i \tilde F_i. \label {equ:5}
\end {eqnarray}
From (\ref {equ:fourth}) and (\ref {equ:5}), we obtain
\bqs
\X_0 \phi_1 &=&\sum _{i=1} ^{n-1} \X_i F_i +\frac 12 G -\sum _{i=1} ^{n-1} \X_i \tilde F_i -\frac 12 \tilde G -\frac 12 b r.
\eqs
Multiplying both sides by the inverse of $\X_0$, we then have \refe{equ:important}. \eop

\begin{description}
\item[Remarks.] Note that $P_1$ is assumed to be symmetric. Therefore quadratically equivalent coefficients $\bar F_i$ ($i=1,\cdots,n-1$) and $\bar G$ have to ensure the symmetry of the right hand side of \refe{equ:important}. Thus \refe{equ:important} constitutes a necessary condition for all equivalent systems of \refe{sys:consin}, including the Brunovsky forms.
\end{description}

\subsection{Computation of the Brunovsky forms and the transformation matrices}
In this section, given $r$ as well as the coefficients of \refe{sys:consin}, $F_i$ ($i=1,\cdots,n$) and $G$, we show how to choose $\bar F_i$ ($i=1,\cdots,n$) and $\bar G$ in Brunovsky forms, which, first, satisfy the necessary condition \refe{equ:important} and, second, have the smallest number of non-zero terms.

Since $\bar F_n$ does not appear in the necessary condition, \refe{equ:important}, we simply let $\bar F_n=0$ in Brunovsky forms. To determine other coefficients, we decompose the terms related to original system \refe{sys:consin} and $r$ on the right hand side of \refe{equ:important} as
\bqn
\X_0^{-1}\left (\sum _{i=1} ^{n-1} \X_i F_i+\frac 12 G-\frac12 b r\right)&=&L+D+U, \label{decompose}
\eqn
where $L,D,U$ are strictly lower triangular, diagonal, and strictly upper triangular matrices respectively. In order to add $\X_0^{-1} (\sum _{i=1} ^{n-1} \X_i \tilde F_i +\frac 12 \tilde G)$ to the right hand side of  \refe{equ:important} to obtain a symmetric matrix, we have, according to the property of the decomposition, the following three cases.

\ben
\item If $L-U^T=0$, $L+D+U$ is already symmetric, and we simply set $\bar F_i$ ($i=1,\cdots,n-1$) and $\bar G$ to be 0. In this case, the Brunovsky form of \refe{sys:consin} is a linear system; i.e., \refe{sys:consin} can be linearized.

\item When $L-U^T\neq 0$, we can set $\X_0^{-1} (\sum _{i=1} ^{n-1} \X_i \tilde F_i +\frac 12 \tilde G)$ to a lower-triangular matrix, $L-U^T$. In this case, $\phi_1=U^T+D+U$. According to the properties of $\X_0$, $\sum _{i=1} ^{n-1} \X_i \tilde F_i +\frac 12 \tilde G$ is a full matrix.
\item When $L-U^T\neq 0$, we can also set $\X_0^{-1} (\sum _{i=1} ^{n-1} \X_i \tilde F_i +\frac 12 \tilde G)$ to be an upper-triangular matrix, $U-L^T$, $\sum _{i=1} ^{n-1} \X_i \tilde F_i +\frac 12 \tilde G$ is a triangular matrix defined by \refe{def:Delta}, which has $n(n-1)/2$ non-zero terms as $U-L^T$. In this case, we have
\bqn
\sum _{i=1} ^{n-1} \X_i \tilde F_i +\frac 12 \tilde G&=&\X_0(U-L^T)\equiv \Delta_1,\label{last}
\eqn
and the first transformation matrix is given by
\bqn
\phi_1&=&L+D+L^T \label{equ:phi1}.
\eqn
\item In addition to the aforementioned cases, $\X_0^{-1} (\sum _{i=1} ^{n-1} \X_i \tilde F_i +\frac 12 \tilde G)$ can be also either $L-U^T$ or $U-L^T$ plus an arbitrary symmetric matrix.
\een

Comparing Cases 2, 3, and 4, we can see that the solutions in Cases 2 and 4 have more non-linear terms than in Case 3. Therefore, in the Brunovsky forms, we set $\X_0^{-1} (\sum _{i=1} ^{n-1} \X_i \tilde F_i +\frac 12 \tilde G)$ to be $U-L^T$, and from \refe{last}, we can obtain two types of Brunovsky forms as follows.

First, by setting $\tilde G=0$, we have from \refe{last} that $\sum _{i=1} ^{n-1} \X_i \tilde F_i =\Delta_1$, which is a lower skew-triangular matrix. To obtain the smallest number of, i.e., $n(n-1)/2$, non-zero items in the Brunovskey forms, we can select $\bar F_i$ to be a diagonal matrix with main diagonal elements as $(0,\cdots,0,\bar f_{i,i+1},\cdots,\bar f_{i,n})$, which can be uniquely computed as follows. From properties of $\X_i$ ($i=1,\cdots,n-1$), we first have \begin{eqnarray}
\bar f_{1,j}=(\X_1 \tilde F_1)_{n-j+2,j}=(\Delta_1)_{n-j+2,j}, j=2,\cdots,n, \label{barf1}
\end{eqnarray} and
\begin{eqnarray}
\Delta_2=\Delta_1-\X_1 \tilde F_1. \label{barf2}
\end{eqnarray}
Then, we can compute $\bar f_{2,j}$ for $j=3,\cdots,n$ from $\Delta_2$ and $\Delta_3=\Delta_2-\X_2 \tilde F_2$ in the same fashion. By repeating this process, we can obtain all non-zero elements in $\tilde F_i$ for $i=1,\cdots,n-1$. Since all these matrices are uniquely determined by $\Delta_1$, the original system is uniquely equivalent to the following system,
\bqn
\dot x_i&=&x_{i+1}+b_i\nu+\sum_{j=i+1}^n \tilde f_{ij} x_j^2 +O(x,\nu)^3\label {norm:c12}, \qquad i=1,\cdots,n,
\eqn
which is type I, complete-quadratic Brunovsky form in \cite{kang1992}.

Second, by setting $\bar F_i=0$ ($i=1,\cdots,n-1$), we have $\tilde G=2\Delta_1$ and the corresponding equivalent system,
\bqn
\dot x_i&=&x_{i+1}+b_i\nu+\sum_{j=n-(i+1)}^n \tilde G_{ij} x_j\nu +O(x,\nu)^3\label {norm:c11}, \qquad i=1,\cdots,n,
\eqn
which is type II, bilinear Brunovsky form in \cite{kang1992}. We can see that $\bar G$ is also  uniquely determined by the original system.

Once we have the Brunovsky forms of a continuous quadratic system, we can solve the corresponding transformation matrices $\phi_i$ ($i=1,\cdots,n$) and $\alpha$ from \refe{equ:phi1}, \refe{equ:third}, and \refe{equ:alpha}. These transformation matrices are also unique corresponding to each Brunovsky form.

\subsection{Discussions}
In the preceding subsection, we finished computing the two types of Brunovsky forms and the corresponding transformation matrices of a linearly controllable quadratic system, \refe{sys:consin}. Besides, we showed that the Brunovsky form of each type and the corresponding transformation matrices are unique. However, the uniqueness is dependent on $r$. That is, for different values of $r$, a system can have multiple solutions of type I or type II Brunovsky forms. A straightforward strategy to prevent the multiplicity is to set $r=0$, and we have a simplified version of transformations, defined by \refe{qud_tran}. I.e., with \refe{qud_tran}, a quadratic system has unique type I and type II Brunovsky forms. Further, one can follow the arguments of Kang and Krener \cite{kang1992} but with $r=0$ \footnote{Due to the difference in notations, $r=0$ is equivalent to saying $\beta^{[1]}=0$ in \cite{kang1992}.} and prove that \refe{qud_tran} still defines quadratic equivalence in the same sense.

\section {Computation of discrete quadratic Brunovsky forms}
In this section, we apply the method proposed in the preceding section and solve Brunovsky forms for a discrete quadratic control system \refe{sys:dissin}, but with simplified transformations defined in \refe{qud_tran}. We carry out our study in the same three-step, constructive procedure as for the continuous system. Compared to continuous system, the discrete system has one more term, and the relationships between coefficients and transformations and the resulted Brunovsky forms are fundamentally different from those of continuous systems.

\subsection{Quadratically equivalent systems of \refe{sys:dissin}}

\begin{definition}\label{def4.1}
Assume $A\in\Re^{n\times n}$ is in the linear Brunovsky form defined in \refe{def:AB}, we define a linear operator $\L:\Re^{n\times n}\to\Re^{n\times n}$ by
\bqn
\ba{lcl}
\L^0 \phi   &=& \phi, \nonumber\\
\L \phi       &=& A^T \phi A,\\
\L^{i+1} \phi &=&\L \cdot \L^i \phi, \quad i=0,1,\cdots.
\ea\label{equ4:3}
\eqn
\end{definition}
\begin{description}
\item[Properties of $\L$.]  The linear operator $\L$ has the following properties:
\ben
\item $\L^i P=0$ when $i\geq n$.
\item The nullity of $\L$ is $2n-1$; if $P\in ker(\L)$, then $P_{ij}= 0$ when $(n-i)(n-j)>0$, and elements in $n$th row and $n$th column can be arbitrarily selected.
\item $\L$ is not invertible.
\item If $P$ is symmetric, $\L P$ is symmetric; but the inverse proposition may not be true.
\een

\item[Remarks.] Note that the linear operator $\L$ for the discrete system is different from that for continuous systems. This fundamental difference leads to the difference in the relationships between coefficients and transformations and the Brunovsky forms.

\end{description}

\begin{theorem}
Discrete system \refe{sys:dissin} is quadratically equivalent, in the sense of \cite{kang1992}, under the quadratic transformations in \refe{qud_tran}, to a system, whose $i$th ($i=1,\cdots,n$) equation is
\bqn
x_i(t+1)&=&x_{i+1}(t)+b_i\nu(t)+x^T(t)\tilde F_i x(t)+\tilde G_i x(t)\nu(t){+O(x,\nu)^3}, \label {norm:ds1}
\eqn
where $x_{n+1}(t)\equiv 0$ is a dummy state variable, $\bar F_1,\cdots,\bar F_n$ are symmetric, $\bar G_i$ is the $i$th row of the matrix $\bar G$, and the coefficients $\bar F_i$ and $h_i$ ($i=1,\cdots,n$) and $\bar G$ are defined by
\bqn
\tilde F_i &=&F_i+\phi_{i+1}-\L\phi_i{-b_i\alpha}, \label {equ4:2}\\
\tilde G_i&=&G_i-2b^T\phi_iA, \label{def:g}\\
\gamma_i&=&(\phi_i)_{nn}, \label {equ4:4}
\eqn
where $P_{n+1}$ is a zero dummy matrix, $b_i$ is the $i$th element of $b$, and $(\phi_i)_{nn}$ is the $(n,n)$ entry of the matrix $\phi_i$.
\end{theorem}

{\em Proof.} The proof is similar to that of Theorem \ref{theorem:con-equiv} and, therefore, omitted. \eop
\begin{description}
\item[Remarks.] Note that \refe{def:g} is different from \refe{equ:second}, and we have an additional equation \refe{equ4:4}.
\end{description}

\subsection{A necessary condition for the equivalent systems}

\begin{lemma}
The mapping from the coefficients $G$ and $\bar G$ to the transformation matrices $P_i$ ($i=1,\cdots,n$) is given by
\begin {eqnarray}
\matphi{n}{n}\cdot A&=&\frac 12 G-\frac 12 \tilde G .\label {equ4:6}
\end {eqnarray}
\end{lemma}

{\em Proof.} The proof is similar to that of Lemma \ref{lemma:pg} and, therefore, omitted. \eop

\begin{lemma}
The mapping from the coefficients $F_i$, $\bar F_i$ ($i=1,\cdots,n$) and $h$ to the transformation matrices $P_i$ ($i=1,\cdots,n$) and $Q$ is given by
\bqn
\phi_{i+1} &=&\L^i\phi_1-\sum_{j=0}^{i-1} \L^j F_{i-j} +\sum_{j=0}^{i-1}\L^j \tilde F_{i-j}, \quad i=1,\cdots, n-1, \label {equ4:11} \\
\alpha  &=&\sum_{j=0}^{n-1} \L^j (F_{n-j}-\tilde F_{n-j}), \label{def:alpha}
\eqn
and
\bqn
\ba{lcl}
(\phi_1)_{n-i\:n-i}&=&\sum_{j=0}^{i-1}(\L^jF_{i-j})_{nn}+\gamma_{i+1}, \qquad i=1,\cdots,n-1,\\
(\phi_1)_{nn}&=&\gamma_1.
\ea\label{eqn:diagonal}
\eqn
\end{lemma}
{\em Proof.} The derivation of \refe{equ4:11} and \refe{def:alpha} from \refe{equ4:2} is similar to that in the proof of Lemma \ref{lemma:pf} and, therefore, omitted. Here we only show how \refe{eqn:diagonal} can be derived as follows. From \refe{equ4:11}, we have ($i=1,\cdots,n-1$)
\bqs
(\phi_{i+1})_{nn}&=&(\L^i\phi_1)_{nn}-\sum_{j=0}^{i-1}(\L^jF_{i-j})_{nn}
        =(\phi_1)_{n-i\:n-i}-\sum_{j=0}^{i-1}(\L^jF_{i-j})_{nn}.
\eqs
Comparing this with \refe{equ4:4}, we then obtain \refe{eqn:diagonal}, from which we can see that the diagonal entries of $P_1$ are uniquely determined by the coefficients $F_i$ ($i=1,\cdots,n$) and $\gamma$. \eop

\begin{definition}
With $A\in\Re^{n\times n}$ and $\L$ given by \refe{def:AB} and Definition \ref{def4.1} respectively, we define a series of linear operators $\X_i:\Re^{n\times n}\to \Re^{n\times n}$ ($i=0,1,\cdots$) by
\bqn
\ba{lcl}
\X_0 \phi &=&\matx nn, \\
\X_i \phi &=& (A^T)^i \X_0 \phi.
\ea
\eqn
\end{definition}
\begin{description}
\item[Properties of $\X_i$.] The linear operators $\X_i$ ($i=0,1,\cdots$) have the following properties:
\ben
\item $\X_i=0$ when $i\geq n$.
\item The rank of $\X_0$ is $n(n+1)/2$; denote the $(i,j)$th entry of $\phi$ by $P_{ij}$ ($i,j=1,\cdots,n$), then the $(i,j)$th entry of $\X_0 P$ can be written as
\bqn
(\X_0 P)_{ij}&=&\cas{{ll} P_{n-i+1,j-i+1},& j\geq i;\\ 0,&j<i.}\label{property}
\eqn
\item From Property 2, $\X_0$ is not invertible. Note that its continuous counterpart in Definition \ref{def:xc} is invertible.
\een
\end{description}

\begin{theorem}
The mapping from the coefficients $F_i$, $G$, $\bar F_i$, and $\bar G$ ($i=1,\cdots,n$) to $P_1$ is given by
\bqn
\X_0 \phi_1 A &=&\sum _{i=1} ^{n-1} \X_i F_i A +\frac 12 G -\sum _{i=1} ^{n-1} \X_i \tilde F_i A -\frac 12 \tilde G\label {equ4:last}
\eqn
\end{theorem}
{\em Proof.} The proof is similar to that of Theorem \ref{con:p1} and omitted. \eop

\begin{description}
\item[Remarks.] According to \refe{property}, we can find the $(i,j)$th entry of $\X_0 P A$,
\bqn
(\X_0 P_1 A)_{ij}&=&\cas{{ll} (P_1)_{n-i+1,j-i},& j> i;\\ 0&j\leq i.} \label{pro:trans}
\eqn
I.e., $\X_0 PA$, the left hand side of \refe{equ4:last}, is an upper-triangular matrix, which contains all the non-diagonal entries of $P_1$. Therefore, $\bar F_i$ ($i=1,\cdots,n$) and $\bar G$ of a quadratically equivalent system of \refe{sys:dissin} must satisfy that the right hand side of \refe{equ4:last} be also upper triangular. Thus, \refe{equ4:last} is a necessary condition for all equivalent systems of \refe{sys:dissin}, including the Brunovsky form.
\end{description}

\subsection{Computation of the Brunovsky form and the transformation matrices}

\begin{theorem}
In the Brunovsky form of \refe{sys:dissin}, the coefficients are
\bqn
\ba{lcl}
\tilde F_i&=&0, \qquad \forall i=1,\cdots,n,\\
\tilde G&=&2 (L+D),
\ea\label{dis:br}
\eqn
where $L$ and $D$ are the strictly lower triangular part and the diagonal part of $\sum _{i=1} ^{n-1} \X_i F_i A +\frac 12 G$.
\end{theorem}

{\em Proof.} The first two terms on the right hand side of \refe{equ4:last} can be written as
\bqn
\sum _{i=1} ^{n-1} \X_i F_i A +\frac 12 G &=& L+D+U, \label{dis:de}
\eqn
where $L,D,U$ are strictly lower triangular, diagonal, and strictly upper triangular matrices respectively. To ensure the right hand side of \refe{equ4:last} to be an upper triangular matrix, we set $\sum _{i=1} ^{n-1} \X_i \tilde F_i A +\frac 12 \tilde G$ to be
\bqn
\sum _{i=1} ^{n-1} \X_i \tilde F_i A +\frac 12 \tilde G =L+D, \label{def:fg}
\eqn
in which there are $n(n+1)/2$ non-zero elements at most, and we hence have
\bqn
\X_0 P_1 A&=& U. \label{dis:p1u}
\eqn

Due to the properties of $\X_i$, no simple $\tilde F_i$'s with $n(n+1)/2$ non-zero elements can satisfy \refe{def:fg}. We thus simply set the coefficients as in \refe{dis:br}, and obtain the only type of Brunovsky form of \refe{sys:dissin}: \bqn
 x_i(t+1)&=&x_{i+1}(t)+b_i\nu(t)+ \sum_{j=1}^i \tilde G_{ij}x_j(t)\nu +O(x,\nu)^3,\quad i=1,\cdots,n.
\eqn
\eop

\begin{description}
\item[Remarks.] The Brunovsky form of discrete systems corresponds to type II form of continuous systems and has $n(n+1)/2$ bilinear terms. After finding the Brunovsky form of a discrete quadratic system, we can solve the corresponding transformation matrices $P_i$ ($i=1,\cdots,n$) and $Q$ as follows. From \refe{pro:trans} and \refe{dis:p1u}, we can find the non-diagonal elements of $P_1$, and the diagonal elements from \refe{eqn:diagonal}. Then $\phi_i$ ($i=2,\cdots,n$) and $\alpha$ can be calculated from \refe{equ4:11} and \refe{def:alpha}.
\end{description}

\section{Algorithms and examples}
In this section, we summarize the algorithms for computing Brunovsky forms and the corresponding transformations of both continuous and discrete linearly controllable quadratic systems with a single input. To prevent multiplicity in Brunovsky forms, we use the simplified version of transformations, \refe{qud_tran}, for both systems. Thus, $r\equiv0$ for all formulas in Section 3. Following each algorithm, an example is given.

\begin{algorithm}\label{alg1}
Computation of continuous quadratic Brunovsky forms.
\bi
\item [Step 1] Compute $\X_0^{-1}\left(\sum _{i=1} ^{n-1} \X_i F_i+\frac 12 G\right)$ and decompose it into $L+D+U$.

\item [Step 2] Compute $\phi_1=L+D+L^T$ and $\Delta_1=\X_0(U-L^T)$.

\item [Step 3] For type I Brunovsky forms, set $\tilde G=0$ and solve the equation
\bqs
\sum _{i=1} ^{n-1} \X_i \tilde F_i=\Delta_1
\eqs
to compute $\tilde F_i$ ($i=1,\cdots,n-1$) as in \refe{barf1} and \refe{barf2}. Go to Step 5.

\item [Step 4] For type II Brunovsky forms, set $\tilde F_i=0$ ($i=1,\cdots,n$) and $\tilde G=2\Delta_1$. Go to  step 5.

\item [Step 5] Compute $\phi_i$ ($i=2,\cdots,n$) and $\alpha$ from \refe{equ:third} and \refe{equ:alpha}.
\ei

\end{algorithm}

\begin{example} Find type I Brunovsky form of a continuous quadratic control system, which is already in type II Brunovsky form,
\bqn
\ba{lcl}
\dot \xi_1&=&\xi_2+O(\xi,\mu)^3,\\
\dot \xi_2&=&\mu+\xi_2 \mu +O(\xi,\mu)^3.
\ea\label{sys1}
\eqn
\end{example}
{\em Solutions.}
In this system, $n=2$, and coefficients are
\bqs
F_1&=&F_2=0,\\
G&=&\mat{{cc}0&0\\0&1}.
\eqs

Following Algorithm \ref{alg1}, we find type I Brunovsky form,
\bqn
\ba{lcl}
\dot x_1&=&x_2 +\frac 12 x_2^2 + O(x,\nu)^3,\\
\dot x_2&=&\nu+O(x,\nu)^3,
\ea\label{sys2}
\eqn
and the corresponding transformations,
\bqn
\ba{lcl}
\xi_1&=&x_1,\\
\xi_2&=&x_2+\frac 12 x_2^2,\\
\mu&=&\nu.
\ea\label{tran1}
\eqn

Reversely, if given \refe{sys2}, we find the transformations for obtaining \refe{sys1},
\bqs
x_1&=&\xi_1,\\
x_2&=&\xi_2-\frac 12 \xi_2^2,\\
\nu&=&\mu,
\eqs
which are the inverses of \refe{tran1}. \eop

\begin{algorithm} \label{alg2}
Computation of discrete quadratic Brunovsky form.
\bi
\item [Step 1] Compute $\sum _{i=1} ^{n-1} \X_i F_i A +\frac 12 G$ and decompose it into $L+D+U$.
\item [Step 2] Compute $\bar G=2(L+D)$.
\item [Step 3] Compute the non-diagonal elements of $\phi_1$ from \refe{pro:trans} and \refe{dis:p1u}.
\item [Step 4] Compute the diagonal elements of $\phi_1$ from \refe{eqn:diagonal}.
\item [Step 5] Compute $\phi_i$ ($i=2,\cdots,n$) and $\alpha$ from from \refe{equ4:11} and \refe{def:alpha}.
\ei
\end{algorithm}

\begin{example} Find the Brunovsky form of the following discrete system:
\bqs
\xi_1(t+1)&=&\xi_2(t)+\xi_1^2(t)+\xi_2^2(t)+\mu^2(t)+O(\xi,\mu)^3,\\
\xi_2(t+1)&=&\mu(t)+\mu^2(t)+O(\xi,\mu)^3.
\eqs
\end{example}
{\em Solutions.} In this system, n=2 and coefficients are
\bqs
F_1&=&\mat{{cc}1&0\\0&1},\\
F_2&=&0,\\
G&=&0,\\
\gamma&=&\mat{{c}1\\1}.
\eqs

Following Algorithm \ref{alg2}, we find the Brunovsky form,
\bqs
x_1(t+1)&=&x_2(t)+O(x,\nu)^3,\\
x_2(t+1)&=&\nu(t)+O(x,\nu)^3,
\eqs
and the corresponding transformations
\bqs
\xi_1(t)&=&x_1(t)+2x_1^2(t)+x_2^2(t),\\
\xi_2(t)&=&x_2(t)-x_1^2(t)+x_2^2(t),\\
\mu(t)&=&\nu(t)-x_2^2(t).
\eqs
Thus, this system is linearized. \eop

\section {Conclusions}
In this paper, we proposed a method for computing Brunovsky forms of both continuous and discrete quadratic control systems, which are linearly controllable with a single input. Our approach is constructive in nature and computes Brunovsky forms explicitly in a three-step manner, with the introduction of linear operators $\L$ and $\X_i$ ($i=0,\cdots,n-1$): (i) we first derived relationships between the coefficients of original systems, those of the quadratically equivalent systems, and the corresponding transformation matrices; (ii) after simplifying these relationships, we obtained a necessary condition for all equivalent systems; (iii) from the necessary conditions, we finally computed the Brunovsky forms and the corresponding transformations. However, the linear operators $\L$ and $\X_i$ and, therefore, the Brunovsky forms are fundamentally different for continuous and discrete systems. For a continuous quadratic system, which has at most $n(n+1)^2/2$ nonlinear terms, there are two types of Brunovsky forms with $n(n-1)/2$ nonlinear terms each; for a discrete quadratic system, there is only one type, and the numbers of nonlinear terms of a general system and its Brunovsky form are $n(n+1)^2/2+n$ and $n(n+1)/2$, respectively.

For continuous systems, we used the same transformations as in \cite{kang1992} and found that the inclusion of transformation vector $r$ yields multiple solutions of a type of Brunovsky forms. Therefore, we suggested to set $r=0$ in order to maintain the uniqueness. I.e., under transformations defined by \refe{qud_tran} for both continuous and discrete systems,  the Brunovsky forms and the corresponding quadratic transformations always exist and are uniquely determined by the original systems. In this sense, our study can be viewed as a constructive proof of the existence and uniqueness of Brunovsky forms.

Finally, the method proposed in this paper could be extended for studying quadratic control systems, not linearly controllable or with multiple input. In addition, the method could be useful for applications involving analysis and control of quadratic systems.

\end {document}